\documentclass[reqno,12pt]{amsart}
\usepackage{graphicx}
\usepackage{amsmath,amssymb}
\usepackage[top=1in,bottom=1in,left=1in,right=1in]{geometry}
\newtheorem{thm}{Theorem}[section]

\newtheorem{lem}[thm]{Lemma}

\theoremstyle{definition}

\theoremstyle{remark}
\newtheorem{rem}{Remark}[section]
\numberwithin{equation}{section}

\begin{document}
\title[Nonlocal parabolic problem]
      {Asymptotic behavior of a nonlocal parabolic problem in Ohmic heating process}%

\author[Liu]{Qilin Liu}%
\address{Department of Mathematics, Shanghai Jiao Tong University,
         Shanghai 200240, PR China}
\email{liuqlseu@yahoo.com.cn}

\author[Liang]{Fei Liang}%
\address{Department of Mathematics, Southeast University, Nanjing, 210096,  Jiangsu, PR China}
\email{liangfei1980@yahoo.com.cn}

\author[Li]{Yuxiang Li}%
\address{Department of Mathematics, Southeast University, Nanjing, 210096,  Jiangsu, PR China}
\email{lieyx@seu.edu.cn}

\thanks{Supported in part by PRC Grants NSFC 10571087, 10671123 and 10601012
        and Southeast University Award Program for Outstanding Young Teachers 2005.}



\begin{abstract}
In this paper, we consider the asymptotic behavior of the nonlocal
parabolic problem
\[
  u_{t}=\Delta u+\displaystyle\frac{\lambda
  f(u)}{\big(\int_{\Omega}f(u)dx\big)^{p}},\ \ \ x\in \Omega,\ t>0,
\]
with homogeneous Dirichlet boundary condition, where $\lambda>0,\
p>0$, $f$ is nonincreasing.  It is found that: (a) For $0<p\leq1$,
$u(x,t)$ is globally bounded and the unique stationary solution is
globally asymptotically stable for any $\lambda>0$; (b) For $1<p<2$,
$u(x,t)$ is globally bounded for any $\lambda>0$; (c) For $p=2$, if
$0<\lambda<2|\partial\Omega|^2$, then $u(x,t)$ is globally bounded,
if $\lambda=2|\partial\Omega|^2$, there is no stationary solution
and $u(x,t)$ is a global solution and $u(x,t)\rightarrow\infty$ as
$t\rightarrow\infty$ for all $x\in\Omega$, if
$\lambda>2|\partial\Omega|^2$, there is no stationary solution and
$u(x,t)$ blows up in finite time for all $x\in\Omega$; (d) For
$p>2$, there exists a $\lambda^*>0$ such that for
$\lambda>\lambda^*$, or for $0<\lambda\leq\lambda^*$ and $u_0(x)$
sufficiently large, $u(x,t)$ blows up in finite time. Moreover, some
formal asymptotic estimates for the behavior of $u(x,t)$ as it blows
up are obtained for $p\geq2$.
\end{abstract}
\maketitle

\section{Introduction}

In this paper, we study the asymptotic behavior of the following
nonlocal parabolic problem
\begin{eqnarray}\label{e:main}
\begin{cases}
  u_{t}=\Delta u+\displaystyle\frac{\lambda
  f(u)}{\big(\int_{\Omega}f(u)dx\big)^{p}},\ \ \ &x\in \Omega,\
  t>0,\\
  u(x,t)=0,            &x\in\partial \Omega,\ t>0,\\
  u(x,0)=u_{0}(x),     &x\in \Omega,
\end{cases}
\end{eqnarray}
where $\lambda>0$ and $u(x,t)=u(x,t,\lambda)$ stands for the
dimensionless temperature of a conductor when an electric current
flows through it \cite{F, L1, T}. $\Omega$ is a bounded domain of
$R^n$ with $C^2$ boundary.  The nonlinear function $f(s)$ satisfies
the condition
\begin{equation}\label{e:f}
  f(s)>0,\ f^\prime(s)<0,\ s\geq0,\ \int_0^{\infty}f(s) ds<\infty,
\end{equation}
and represents, depending on the problem, either the electrical
conductivity or the electrical resistance of the conductor.
Condition (\ref{e:f}) permits us to use comparison methods, see
\cite{L1, L2, S, T}.  Also for simplicity, we assume $u_0(x)$ is
continuous with $u_0(x)=0,x\in\partial\Omega$ and
$u_0(x)\geq0,x\in\Omega$. Without loss of generality, we may assume
that $\int_0^{\infty}f(s) ds=1$.

A similar study had been undertaken in \cite{L1, L2, T, K, KT} for
the nonlocal reaction-diffusion problem (\ref{e:main}) for $p=2$.
Lacey \cite{L1, L2} and Tzanetis \cite{T} proved the occurrence of
blow-up for the one-dimensional problem and for the two-dimensional
radially symmetric problem, respectively. First they estimated the
supremum $\lambda^{*}$ of the spectrum of the related steady-state
problem and then they proved the blow-up, for $\lambda>\lambda^{*}$,
by constructing some blowing-up lower solutions.  Using some ideas
of Bebernes and Lacey \cite{BL}, Kavallaris and Tzanetis \cite{KT}
generalized the blow-up results for $\lambda>\lambda^*$ and
dimensions $n\geq2$ if $u_0$ is sufficiently large and $f(s)$
satisfies
\begin{equation}\label{e:condition}
  \int_0^\infty[sf(s)-s^2f'(s)]ds<\infty.
\end{equation}
Kavallaris and Lacey \cite{K} showed that the solution
$u^*(x,t)=u(x,t,\lambda^*)$ is global in time and diverges in the
sense $\mid\mid u^*(\cdot,t)\mid\mid_\infty\rightarrow\infty$ as
$t\rightarrow\infty$ when $n=1, \Omega=(-1,1)$ and $f(s)$ satisfies
(\ref{e:f}) or $n=2, \Omega=\{(x,y)\in R^2:x^2+y^2<1\}$ and
$f(s)=e^{-s}$. Moreover, it is proved that this divergence is
global, i.e. $u^*(x,t)\rightarrow\infty$ as $t\rightarrow\infty$ for
all $x\in\Omega$.

Throughout this paper, we always assume that the domain $\Omega$
satisfies the following condition:

(H) $\Omega\in R^n$ is a convex smooth bounded domain; for any point
$y_0\in\partial\Omega$, there exists a hyperplane $S_{y_0}$ such
that $S_{y_0}$ is tangent to $\Omega$ at $y_0$
($\{y_0\}=S_{y_0}\cap\partial\Omega$).

The main purpose of this paper is to generalize and improve the
results for dimensions $n\geq2$ and $p>0$ obtained in \cite{L1, L2,
T, K, KT}. Our main results read as follows.
\begin{enumerate}
  \item[$\bullet$] If $0<p\leq1$, then $u(x,t)$ is globally bounded and there exists a unique stationary solution which is
                   globally asymptotically stable for any $\lambda>0$.
  \item[$\bullet$] If $1<p<2$, then $u(x,t)$ is globally bounded for any $\lambda>0$.
  \item[$\bullet$] Assume $p=2$, let $\lambda^*=2|\partial\Omega|^2$.  If $0<\lambda<\lambda^*$, $u(x,t)$ is globally bounded.  If
                   $\lambda=\lambda^*$, there is no stationary solution
                   and $u^*(x,t)$ is a global-in-time solution and
                   $u^*(x,t)\rightarrow\infty$ as $t\rightarrow\infty$ for all
                   $x\in\Omega$.  If $\lambda>\lambda^*$, there is
                    no stationary solution and $u(x,t)$ blows up globally in finite time $T$ but the condition (\ref{e:condition}) and $u_0$ sufficiently
                   large are not required.
  \item[$\bullet$] If $p>2$, then
                   there exists a critical value $\lambda^*$ such that for
                   $\lambda>\lambda^*$ or for any $0<\lambda\leq\lambda^*$ and $u_0(x)$ sufficiently large, $u(x,t)$ blows up globally in finite time $T$.
  \item[$\bullet$] We also obtain some formal asymptotic estimates for the local
                   behavior of $u(x,t)$ as it blows up for $p\geq2$.
\end{enumerate}

This paper is organized as follows.  In Section $2$ we consider the
steady-state problem corresponding to (\ref{e:main}).  In Section
$3$, we investigate the behavior of some critical solutions of the
equation (\ref{e:main}) for $p=2$.  Section $4$ is devoted to some
formal asymptotic estimates for the local behavior of $u(x,t)$ as it
blows up in finite time for $p\geq2$.

\section{Steady-state problem}

The steady states of the problem (\ref{e:main}) play an important
role in the description of the asymptotic behavior of the solutions
of (\ref{e:main}) and the construction of the lower and upper
solutions, so we first consider the stationary problem of
(\ref{e:main}).  The stationary problem corresponding to
(\ref{e:main}) is
\begin{equation}\label{e:steady1}
  \Delta w+\frac{\lambda
  f(w)}{\big(\int_{\Omega}f(w)dx\big)^{p}}=0,\ \ \ x\in\Omega;\ \ \ w=0,\ \ \ x\in\partial\Omega.
\end{equation}

In order to study the nonlocal problem (\ref{e:steady1}), let us
first consider the following local problem:
\begin{equation}\label{e:steady2}
  \Delta w+\mu f(w)=0,\ \ \ x\in\Omega;\ \ \ w=0,\ \ \
  x\in\partial\Omega,
\end{equation}
where $\mu\geq0$ and $f(s)$ satisfies (\ref{e:f}).  It is well-known
that the basic theory of monotone schemes can be carried out for the
problem (\ref{e:steady2}).  Therefore, there exists a solution in
$H_0^1(\Omega)$.  Moreover, the straightforward argument, based on
the coercivity of $-\Delta$ with Dirichlet boundary condition,
implies that (\ref{e:steady2}) has a unique positive solution
$w_{\mu}^\Omega$ in $H_0^1(\Omega)$.  The above arguments are
classical and known in the literature \cite{GT}.

In order to establish a relationship between the local problem
(\ref{e:steady2}) and the nonlocal problem (\ref{e:steady1}), we
define a real function $\lambda(\mu)$ by
\begin{equation}\label{e:parameter}
  \lambda(\mu)=\mu\big(\int_\Omega f(w_{\mu}^\Omega)dx\big)^p,
\end{equation}
for any $\mu\geq0$.  This function is well defined due to the
positive character of $w_{\mu}^\Omega$.  From the analyticity of the
solutions $w_{\mu}^\Omega$ on $\mu$, we deduce that the function
$\lambda(\mu)$ is analytical on $\mu$.

It is easy to see the relation between the solutions of problem
(\ref{e:steady2}) and the problem (\ref{e:steady1}).
\begin{thm}\label{lem:2}
If $w$ is a solution of problem (\ref{e:steady1}) for
$\lambda=\lambda_0$, then $w$ is a solution of problem
(\ref{e:steady2}) for $\mu=\lambda_0\big/\big(\int_\Omega
f(w)dx\big)^p$.  Conversely, if $w$ is a solution of problem
(\ref{e:steady2}) for $\mu=\mu_0$, then $w$ is a solution of problem
(\ref{e:steady1}) for $\lambda=\lambda(\mu_0)$.
\end{thm}

Theorem \ref{lem:2} allows us to study problem (\ref{e:steady1}) by
analyzing the behavior of the function $\lambda(\mu)$.  This is the
key idea to solve problem (\ref{e:steady1}).  Now we give some
qualitative properties of the profile of the bifurcation diagram of
the local problem (\ref{e:steady2}).

\begin{lem}\label{prop:1}
Let $w_{\mu}^\Omega$ be the solution of (\ref{e:steady2}), then
\begin{enumerate}
  \item $\partial
        w_{\mu}^\Omega/\partial\mu>0$ for $x\in\Omega$.
  \item  $\lim_{\mu\rightarrow\infty}w_{\mu}^\Omega(x)/\Phi_1^\Omega(x)\rightarrow\infty$,
uniformly in $\Omega$, where $\Phi_1^\Omega(x)$ is the first
normalized eigenfunction of $-\Delta$ in $H_0^1(\Omega)$.
\end{enumerate}
\end{lem}

The proof follows the same line as in \cite{C}, so we omit it.

Now we are going to  prove that the solution of (\ref{e:steady1}) is
unique for any $0<p\leq1$.
\begin{thm}\label{prop:2}
For any $0<p\leq1$, there exists a unique solution of the problem
(\ref{e:steady1}) for any $\lambda\geq0$.
\end{thm}
\begin{proof}
Let us prove that $\lambda(\mu)$ is strictly increasing. Integrating
the equation (\ref{e:steady2}) over $\Omega$, we have
\[
  \int_{\partial\Omega}\frac{\partial
  w}{\partial\nu}ds+\lambda^{\frac{1}{p}}\mu^{\frac{p-1}{p}}=0,
\]
where $\partial/\partial \nu$ is the outward normal derivative,
which implies
\begin{equation}\label{e:parameter1}
  \lambda(\mu)=\mu^{1-p}(-\int_{\partial\Omega}\frac{\partial
  w}{\partial\nu}ds)^p.
\end{equation}
By $0<p\leq1$, $w_\mu=0$ on $\partial\Omega$ and Lemma~\ref{prop:1},
we get
\[
 \lambda'(\mu)>0 \ {\rm for} \  \mu>0\ {\rm and}\
 \lim_{\mu\rightarrow\infty}\lambda(\mu)=\infty.
\]
The proof is completed.
\end{proof}

The following results give us a way to construct sub-solution of
$w_{\mu}^\Omega$ in order to estimate from above the function
$\lambda(\mu)$.
\begin{lem}\label{prop:3}
Let $\Omega'\subset\Omega$.  Then $w_{\mu}^{\Omega'}\leq
w_{\mu}^{\Omega}$ on $\Omega'$ for any $\mu>0$.
\end{lem}

We omit the proof.

We need a lemma concerning the solution to the problem on a ball
\begin{equation}\label{e:steady3}
  \Delta w+\mu f(w)=0,\ \ \ x\in B;\ \ \ w=0,\ \ \
  x\in\partial B.
\end{equation}

\begin{lem}\label{prop:4}
(See \cite[Lemma 5.1]{T})  Let $f(s)$ satisfy (\ref{e:f}),
$\int_0^\infty f(s)ds=1$ and $w_{\mu}^B$ is a solution of
(\ref{e:steady3}), then we have
\begin{equation}\label{e:auxiliary}
 -\frac{1}{\sqrt{\mu}}\frac{d w_{\mu}^B(R)}{dr}<\sqrt{2},\ \ \
 -\lim_{\mu\rightarrow\infty}\frac{1}{\sqrt{\mu}}\frac{d
 w_{\mu}^B(R)}{dr}=\sqrt{2},
\end{equation}
\end{lem}
where $B=\{x\in R^n:\mid x-x_0\mid<R\}$, $r=\mid x-x_0\mid$.

\begin{thm}\label{theo:steady1}
Let $f(s)$ satisfy (\ref{e:f}), $\int_0^\infty f(s)ds=1$, and
$\Omega$ is a bounded domain satisfying (H). Then the following
assertions hold.
\begin{enumerate}
  \item For $1<p<2$, there exists at least one solution of the
        problem (\ref{e:steady1}) for any value $\lambda>0$.
  \item For $p=2$, let $\lambda^*=2|\partial\Omega|^2$, then there exists
        at least one solution of the problem (\ref{e:steady1})
        for $0<\lambda<\lambda^*$ and no solution for
        $\lambda\geq\lambda^*$.  Moreover,
        $\lambda(\mu)<2|\partial\Omega|^2$ for $\mu>0$ and
        $\lim_{\mu\rightarrow\infty}\lambda(\mu)=2|\partial\Omega|^2$.
  \item For $p>2$, there exists a critical value
        $\lambda^*>0$ such that there exist
        at least two solutions of the problem (\ref{e:steady1})
        for $0<\lambda<\lambda^*$, at least one solution for
        $\lambda=\lambda^*$ and no solution for
        $\lambda>\lambda^*$.  Moreover,
        $\lim_{\mu\rightarrow\infty}\lambda(\mu)=0$.
\end{enumerate}
\end{thm}
\begin{proof}
Let $y_0\in\partial\Omega$. Without loss of generality we assume
that $y_0=0$ and the hyperplane $\{x\in R^n:x_1=0\}$ is tangent to
$\Omega$ at $y_0$. By (H), there exist two balls
$\Omega_1,\Omega_2(\Omega_1\subset\Omega\subset\Omega_2)$ which are
tangent to $\Omega$ at $y_0$, where $\Omega_i=\{x\in
R^n:|x-y_i|<R_i,y_i=(L_i,0')\}$. Lemma~\ref{prop:3} implies that
$w_{\mu}^\Omega\geq w_{\mu}^{\Omega_1}$ on $\Omega_1$ and
$w_{\mu}^{\Omega_2}\geq w_{\mu}^{\Omega}$ on $\Omega$. Applying
Lemma~\ref{prop:4}, we
 conclude that
\[
 \sqrt{2}>-\frac{1}{\sqrt{\mu}}\frac{d
 w_{\mu}^{\Omega_2}(0)}{dx_1}\geq-\frac{1}{\sqrt{\mu}}\frac{d
 w_{\mu}^{\Omega}(0)}{dx_1}\geq-\frac{1}{\sqrt{\mu}}\frac{d
 w_{\mu}^{\Omega_1}(0)}{dx_1}, \ \ \mu>0
\]
and
\[
 \sqrt{2}=-\lim_{\mu\rightarrow\infty}\frac{1}{\sqrt{\mu}}\frac{d
 w_{\mu}^{\Omega_2}(0)}{dx_1}\geq-\lim_{\mu\rightarrow\infty}\frac{1}{\sqrt{\mu}}\frac{d
 w_{\mu}^{\Omega}(0)}{dx_1}\geq-\lim_{\mu\rightarrow\infty}\frac{1}{\sqrt{\mu}}\frac{d
 w_{\mu}^{\Omega_1}(0)}{dx_1}=\sqrt{2},
\]
which imply
\[
  -\frac{1}{\sqrt{\mu}}\frac{d
 w_{\mu}^{\Omega}(0)}{dx_1}< \sqrt{2}\ \ \  {\rm
 and}\ \ -\lim_{\mu\rightarrow\infty}\frac{1}{\sqrt{\mu}}\frac{d
 w_{\mu}^{\Omega}(0)}{dx_1}=\sqrt{2}.
\]
Since $y_0$ is arbitrary, it follows that
\[
  -\frac{1}{\sqrt{\mu}}\int_{\partial\Omega}\frac{\partial
  w_{\mu}^\Omega}{\partial\nu}ds<\sqrt{2}|\partial\Omega|\
  \ \ {\rm for}\ \ \mu>0\ \  {\rm and}\ \ -\lim_{\mu\rightarrow\infty}\frac{1}{\sqrt{\mu}}\int_{\partial\Omega}\frac{\partial
  w_{\mu}^\Omega}{\partial\nu}ds=\sqrt{2}|\partial\Omega|.
\]
By (\ref{e:parameter1}), we obtain
\begin{enumerate}
  \item[(i)] If $0<p<2$, then
  $\lim_{\mu\rightarrow\infty}\lambda(\mu)=\infty$.
  \item[(ii)] If $p=2$, then $\lambda(\mu)<2|\partial\Omega|^2$ for
           $\mu>0$ and
           $\lim_{\mu\rightarrow\infty}\lambda(\mu)=2|\partial\Omega|^2$.
  \item[(iii)] If $p>2$, then $\lim_{\mu\rightarrow\infty}\lambda(\mu)=0$.
\end{enumerate}
The proof is completed.
\end{proof}

Let $\mu$ in (\ref{e:steady2}) be a function of $t$. Now we give
some conditions of $\mu(t)$ in order for $w(x;\mu(t))$ to be a lower
or an upper solution of (\ref{e:main}).  We first give a lemma.

\begin{lem}\label{lem:upper-lower}
$w(x;\mu)$ is the solution of (\ref{e:steady2}), then $w_\mu>0$ in
$\Omega$ and $w_\mu$ is bounded.
\end{lem}

\begin{proof}
$w_\mu$ satisfies
\begin{eqnarray*}
\begin{cases}
  \Delta w_\mu+ \mu f'(w)w_{\mu}+f(w)=0
  , \ \ \ &x\in \Omega,\\
  w_\mu(x)=0,                   \ \ \ &x\in\partial \Omega.
\end{cases}
\end{eqnarray*}
Since $f'(s)<0$, the coefficient of $w_\mu$ in this equation is
negative. By the maximum principle we obtain $w_\mu>0$.  Also
$w_\mu$ is finite, indeed for a fixed $\mu$, any sufficiently large
constant is an upper solution, $0\leq w_\mu\leq C$.
\end{proof}

Using Lemma~\ref{lem:upper-lower},
\[
  \inf_{x\in\Omega}\frac{f(w)}{w_\mu}>0,
\]
since $f(w)$ is bounded and away from zero.  Denote $v(x,t)=
w(x;\mu(t))$, then
\begin{equation*}
\begin{split}
  v_{t}-\Delta v-\frac{\lambda f(v)}{\big(\int_\Omega f(v)dx\big)^{p}}
  =w_{\mu}\mu'(t)-\frac{\big(\lambda-\lambda(\mu)\big)f(w)}{\big(\int_\Omega f(w)dx\big)^{p}}.
\end{split}
\end{equation*}
Let $\mu(t)$ be the solution of
\begin{equation}\label{e:upper solution }
  \mu'(t)=\frac{\lambda-\lambda(\mu)}{\big(\int_\Omega f(w)dx\big)^{p}}\inf_{x\in\Omega}\frac{f(w)}{w_\mu},\ \ \ \mu(0)=\mu_0.
\end{equation}

If there exists $\mu_0$ such that
\[
  \lambda\leq\lambda(\mu_0)\ \ {\rm and} \ \ w(x;\mu_0)\geq
  u_{0}(x),
\]
then $v(x,t)$ is decreasing and satisfies
\[
  v_{t}-\Delta v-\frac{\lambda f(v)}{\big(\int_\Omega f(v)dx\big)^{p}}\geq0.
\]
So $v(x,t)$ is a decreasing upper solution of (\ref{e:main}).

If there exists $\mu_0$ such that
\[
  \lambda\geq\lambda(\mu_0)\ \ {\rm and} \ \ w(x;\mu_0)\leq
  u_{0}(x),
 \]
then $v(x,t)$ is increasing and satisfies
\[
  v_{t}-\Delta v-\frac{\lambda f(v)}{\big(\int_\Omega f(v)dx\big)^{p}}\leq0.
\]
So $v(x,t)$ is an increasing lower solution of (\ref{e:main}).

The above preparations in hand, we can discuss the behavior of the
solution of (\ref{e:main}).

\begin{thm}
Assume $0<p\leq1$, then the solution of $(\ref{e:main})$ is globally
bounded and the unique steady state is globally asymptotically
stable for any $\lambda>0$.
\end{thm}
\begin{proof}
From Theorem~\ref{prop:2}, for fixed $\lambda$, there is a unique
steady state $w(x;\mu_1)$ of (\ref{e:steady1}) with
$\lambda=\lambda(\mu_1)$. Take $\overline{\mu}(t)$ satisfying
(\ref{e:upper solution }) with $\mu(0)=\overline{\mu}_0$. For any
initial data $u_{0}(x)>0$, we can select $\overline{\mu}_0$ to
satisfy $w(x;\overline{\mu}_0)\geq u_0(x)$.  This can clearly be
done if we require that $u_0(x)$ and $u'_0(x)$ are
bounded(see\cite{L2}). We also choose $\overline{\mu}_0>\mu_1$.
Since $\lambda'(\mu)>0$, we have
\[
  \lambda\leq\lambda(\overline{\mu}_0)\ \ {\rm and} \ \ w(x;\overline{\mu}_{0})\geq
  u_{0}(x),
\]
thus $\overline{\mu}(t)$ is decreasing and
$\overline{\mu}(t)\rightarrow\mu_1$ as $t\rightarrow\infty$.
 So
$\overline{v}(x,t)=w(x;\overline{\mu}(t))$ is a decreasing upper
solution of the problem and
\[
  \overline{v}(x,t)\rightarrow w(x;\mu_{1}),\ \ \ {\rm as}\ \
  t\rightarrow+\infty.
\]

On the other hand, take $\underline{\mu}(t)$ satisfying
(\ref{e:upper solution }) with $\mu(0)=\underline{\mu}_0$. Since
\[
  \lambda(\mu)\rightarrow 0, \ \ \ {\rm as}\ \
  \mu\rightarrow 0,
\]
we can select $\underline{\mu}_0$ sufficiently small such that
\[
   \lambda\geq\lambda(\underline{\mu}_0)\ \ {\rm and} \ \
w(x;\underline{\mu}_{0})\leq
   u_{0}(x),
\]
thus $\underline{\mu}(t)$ is increasing and
$\underline{\mu}(t)\rightarrow\mu_1$ as $t\rightarrow\infty$. So
$\underline{v}(x,t)=w(x;\underline{\mu}(t))$ is an increasing lower
solution of the problem and
\[
  \underline{v}(x,t)\rightarrow w(x;\mu_{1}),\ \ \ {\rm as}\ \
  t\rightarrow+\infty.
\]

Since $\underline{v}\leq u(x,t)\leq\overline{v}$ and both
$\underline{v}$ and $\overline{v}$ tend to $w(x;\mu_{1})$ as
$t\rightarrow\infty$, we see that $u(x,t)$ exists globally and
$u(x,t)\rightarrow w(x;\mu_{1})$ as $t\rightarrow\infty$.  The above
procedure holds for any initial data $u_{0}(x)$, from which it
follows that the solution $w(x;\mu_{1})$ is globally asymptotically
stable.  The proof is completed.
\end{proof}

\begin{thm}\label{th:steady2}
If $1<p<2$ and $\int_0^\infty f(s)ds=1$, then $u(x,t)$ is globally
bounded for any $\lambda>0$.
\end{thm}
\begin{proof}
For the global boundedness of $u(x,t)$, it suffices to construct an
upper solution which is globally bounded. Select $\mu_0$ so large
that
\[
  \lambda\leq\lambda(\mu_0)\ \ {\rm and} \ \
w(x;\mu_0)\geq
  u_{0}(x),
\]
then $\mu(t)$, the solution of (\ref{e:upper solution }), is
decreasing and therefore $w\big(x;\mu(t)\big)$ is a globally bounded
upper solution.
\end{proof}

For $p=2$, we have the similar result.
\begin{thm}
If $p=2$, $\int_0^\infty f(s)ds=1$ and
$0<\lambda<2|\partial\Omega|^2$, then $u(x,t)$ is globally bounded
for any initial data.
\end{thm}

\section{Behavior of solutions of problem
(\ref{e:main}) for $p=2$}

In this Section, we study the behavior of solutions of the following
nonlocal parabolic problem:
\begin{eqnarray}\label{e:critical}
\begin{cases}
  u_{t}=\Delta u+\displaystyle\frac{2|\partial\Omega|^2
  f(u)}{\big(\int_{\Omega}f(u)dx\big)^{2}},\ \ \ &x\in \Omega,\
  t>0,\\
  u(x,t)=0,            &x\in\partial \Omega,\ t>0,\\
  u(x,0)=u_{0}(x),     &x\in \Omega,
\end{cases}
\end{eqnarray}
where $f$ satisfies (\ref{e:f}) and $\int_0^\infty f(s)ds=1$. By
Theorem~\ref{theo:steady1}, it follows that
$\lambda(\mu)<2|\partial\Omega|^2$ for all $\mu>0$, then we can find
an increasing lower solution $v=w(x;\mu(t))$ with $\mu(t)
\rightarrow\infty$ as $t\rightarrow T\leq\infty$. Thus $u(x,t)$ is
unbounded.  Moreover, $u(x,t)$ is globally unbounded.  Indeed, if
$T=\infty$, from Lemma~\ref{prop:1}, $u(x,t)$ is globally unbounded;
if $T<\infty$, $u(x,t)$ is globally blow-up(see the proof of
Theorem~\ref{theo:blowup} for details).

Now we will prove that $||u(\cdot,t)||_\infty\rightarrow\infty$ as
$t\rightarrow\infty$, i.e. $T=\infty$.  It is sufficient to
construct an upper solution $V(x,t)$ to problem (\ref{e:critical})
which is global in time and unbounded. Without loss of generality,
we assume that the hyperplane $\{x:x_1=1\}$ is tangent to $\Omega$
at $(1,0')$, and $\Omega$ lies in the half-space $\{x:x_1<1\}$. Let
$d(x)=\mathrm{dist}(x,\partial\Omega)$. Set
\begin{eqnarray}\label{e:upper}
\begin{cases}
  V(x,t)=w\big(y(x,t);\mu(t)\big),\ \ 0\leq
  d(x)\leq\varepsilon(t),\ \ x\in\Omega,t>0\\
  V(x,t)=M(t)=\max_{0\leq
  d(x)\leq\varepsilon(t)}w\big(y(x,t);\mu(t)\big),\  \
  d(x)\geq\varepsilon(t), \ \ x\in\Omega,t>0,
\end{cases}
\end{eqnarray}
where $0\leq y(x,t)=d(x)/\varepsilon(t)\leq1$, $\varepsilon(t)>0$ is
a function to be chosen later and $w\big(y(x,t);\mu(t)\big)$
satisfies
\begin{equation}\label{e:upper1}
  w_{yy}+\mu(t)f(w)=0,\ \ \ 0<y<1,\ t>0;\ \ \ w(0;\mu(t))=w'(1;\mu(t))=0,
\end{equation}
or equivalently
\begin{equation}\label{e:upper2}
  w_{rr}+\frac{\mu(t)}{\varepsilon^2(t)}f(w)=0,\ \ \ r=d(x),\ \ \
  0\leq r\leq \varepsilon(t), \ t>0;\ \ \
  w(0)=\frac{dw}{dr}|_{r=\varepsilon(t)}=0,
\end{equation}
and
\begin{eqnarray}\label{e:upper3}
\begin{cases}
  \displaystyle\Delta w-\frac{\Delta
    d}{\varepsilon}\frac{dw}{dy}+\frac{\mu}{\varepsilon^2}f(w)=0,\  0\leq d(x)\leq\varepsilon(t), \
    t>0;\\
   w(y(x,t);\mu(t))=0,\ x\in\partial\Omega,\ t>0,
    \ \frac{dw}{dr}|_{r=\varepsilon(t)}=0,
\end{cases}
\end{eqnarray}
which implies
\begin{eqnarray}\label{e:upper4}
\begin{cases}
  \displaystyle\frac{d^2w\big(y(x_1,0');\mu(t)\big)}{dx_1^2}+\frac{\mu}{\varepsilon^2}f\big(w(y(x_1,0');\mu(t))\big)=0,
    \delta(t)<x_1<1,\ t>0;\\
  \displaystyle w\big(y(1,0');\mu(t)\big)=0,\ \ \
    \frac{dw\big(y(\delta(t),0');\mu(t)\big)}{dx_1}=0,
\end{cases}
\end{eqnarray}
where $\varepsilon(t)=1-\delta(t)$.

From the definition of $w$, it is obvious that $w$, $w_r$ are
continuous at $r=\varepsilon(t)$. We can choose $\mu(0)$(or
equivalently $M(0)$) sufficiently large so that $V(x,0)\geq
u_0(x)$(such a choice is possible since $w\rightarrow\infty$ as
$\mu\rightarrow\infty$ and provided that $u_0(x),u'_0(x)$ are
bounded).

For any $\varepsilon>0$, set
$\Omega_\varepsilon=\{x\in\Omega:0<d(x)<\varepsilon(t)\}$. To prove
that $V(x,t)$ is an upper solution, we need some preliminary
results.

Problem (\ref{e:upper2}) and (\ref{e:upper4}) imply that
\begin{equation}\label{e:upper5}
  w_r(0)=\frac{\sqrt{2\mu}}{\varepsilon}\sqrt{\int_0^M f(s)ds},
\end{equation}
and
\begin{equation}\label{e:upper6}
  \int_{\delta(t)}^1f\big(w(y(x_1,0');\mu(t))\big)dx_1=-\frac{\varepsilon^2}{\mu}\frac{dw\big(y(1,0');\mu(t)\big)}{dx_1}.
\end{equation}
From (\ref{e:upper2}), we get
\begin{equation}\label{e:upper7}
  \frac{w_r}{\sqrt{F(w)-F(M)}}=\frac{\sqrt{2\mu}}{\varepsilon},
\end{equation}
where $F(s)=\int_s^\infty f(\sigma)d\sigma>0$.  Relation
(\ref{e:upper7}) gives
\begin{equation}\label{e:upper8}
  \sqrt{\mu(M)}=\frac{\sqrt{2}}{2}\int_0^M\frac{ds}{\sqrt{F(s)-F(M)}}.
\end{equation}
For $s\leq M$, we have $F(s)-F(M)=f(\theta)(M-s),\ \theta\in[s,M]$
and due to $f'(s)<0$ for $s\geq0$, we get
\begin{equation}\label{e:upper00}
  (M-s)f(M)\leq F(s)-F(M)\leq (M-s)f(s).
\end{equation}
Then
\[
  \sqrt{\mu(M)}\leq
  \frac{\sqrt{2}}{2}\int_0^M(M-s)^{-\frac{1}{2}}f^{-\frac{1}{2}}(M)ds\leq\sqrt{\frac{2M}{f(M)}},
\]
and hence
\begin{equation}\label{e:upper9}
  \mu(M)f(M)\leq 2M\ \ \ {\rm for}\ \ M>0.
\end{equation}
However,
\[
  Mf(M)\leq 2\int_{M/2}^M f(s)ds\leq 2\int_{M/2}^\infty
  f(s)ds\ \ \ {\rm and}\ \ \int_{M/2}^\infty  f(s)ds\rightarrow0\ \
  \ {\rm as}\ \ M\rightarrow\infty,
\]
so $Mf(M)\rightarrow0$ as $M\rightarrow\infty$ and due to
(\ref{e:upper9}) we finally get
\begin{equation}\label{e:upper10}
  \sqrt{\mu(M)}f(M)\rightarrow0\ \ \ {\rm as}\ \ M\rightarrow\infty.
\end{equation}

Next we claim that
$\lim_{\mu\rightarrow\infty}\sqrt{2\mu}/M=\infty$.  Indeed, by
(\ref{e:f}) and (\ref{e:upper8}), we obtain
\[
  \frac{\sqrt{2\mu}}{M}\geq
  \frac{\int_0^M(M-s)^{-\frac{1}{2}}f^{-\frac{1}{2}}(s)ds}{M}=\int_0^1\frac{s^{\frac{1}{2}}(1-s)^{-\frac{1}{2}}}{(Ms
  f(Ms))^{\frac{1}{2}}}ds.
\]
Taking into account $sf(s)\rightarrow0$ as $s\rightarrow\infty$, we
deduce that $\lim_{\mu\rightarrow\infty}\sqrt{2\mu}/M=\infty$, i.e.
\begin{equation}\label{e:upper11}
  \lim_{M\rightarrow\infty}M/\sqrt{2\mu}=0.
\end{equation}

As is indicated in \cite{BL}, $d(x)$ is smooth and more precisely
$|\Delta d|\leq K$, for some $K$, in a neighborhood of the boundary
if $\partial\Omega$ is smooth.  In particular, such a neighborhood
$\Omega_\varepsilon$ consists of all $x\in\Omega$ such that
$d(x,\partial\Omega)\leq\varepsilon(t)$ where $\varepsilon(t)$ is
chosen small enough.

Integrating (\ref{e:upper3}) over $\Omega_\varepsilon$ we obtain
\begin{equation*}
\begin{split}
  \int_{\Omega_\varepsilon}f(w)dx&=-\frac{\varepsilon^2}{\mu}\int_{\partial\Omega}\frac{\partial
  w}{\partial\nu}ds+\frac{\varepsilon}{\mu}\int_{\Omega_\varepsilon}\Delta
  d \frac{dw}{dy} dx\\
  &=\frac{\varepsilon^2|\partial\Omega|}{\mu}w_r(0)+\frac{\varepsilon}{\mu}\int_{\Omega_\varepsilon} \Delta
  d\frac{dw}{dy}dx\\
  &=\varepsilon|\partial\Omega|\sqrt{\frac{2}{\mu}}\sqrt{\int_0^M
  f(s)ds}+\frac{\varepsilon}{\mu}\int_{\Omega_\varepsilon}\Delta
  d \frac{dw}{dy} dx\\
  &\geq \varepsilon|\partial\Omega|\sqrt{\frac{2}{\mu}}\sqrt{\int_0^M
  f(s)ds}-\frac{\varepsilon
  K}{\mu}\int_{\Omega_\varepsilon}\frac{dw}{dy}dx\ \ ({\rm Using}\ \
  \frac{dw}{dy}\geq0)\\
  &\geq \varepsilon|\partial\Omega|\sqrt{\frac{2}{\mu}}\sqrt{\int_0^M
  f(s)ds}+\frac{\varepsilon^2|\partial\Omega|
  K}{\mu}\int_{\delta(t)}^{1}\frac{dw\big((x_1,0');\mu(t)\big)}{dx_1}dx_1\\
  &= \varepsilon|\partial\Omega|\sqrt{\frac{2}{\mu}}\sqrt{\int_0^M
  f(s)ds}-\varepsilon^2|\partial\Omega|
  K\frac{M}{\mu},
\end{split}
\end{equation*}
which implies
\begin{equation}\label{e:upper12}
\begin{split}
  \int_\Omega
  f(V)dx&=\int_{\Omega\setminus\Omega_\varepsilon}f(M)dx+\int_{\Omega_\varepsilon}f(w)dx\\
  &\geq |\Omega\setminus\Omega_\varepsilon|f(M)+\varepsilon|\partial\Omega|\sqrt{\frac{2}{\mu}}\sqrt{\int_0^M
  f(s)ds}-\varepsilon^2|\partial\Omega|
  K\frac{M}{\mu}.
\end{split}
\end{equation}

Our construction of upper solution $V$ depends strongly on the
behavior of the function
\[
  g(s)=\frac{f(s)\sqrt{\mu(s)}}{F(s)}>0.
\]
Since (\ref{e:upper10}) holds and $F(M)\rightarrow0$ as
$M\rightarrow\infty$, we distinguish two cases for the behavior of
$g(M)$.  More precisely the following holds:

\begin{thm}\label{theo:critical}
Let $f(s)$ satisfy (\ref{e:f}), $\int_0^\infty f(s)ds=1$,
$\liminf_{s\rightarrow\infty}g(s)>C>0$ and
$\lim_{s\rightarrow\infty}\mu(s)f(s)=C_0>0$\ (e.g. $f(s)=e^{-s}$) or
$\liminf_{s\rightarrow\infty}\mu(s)f(s)/s=C_1>0$\ ($C_1\leq2$, e.g.
$f(s)=b(1+s)^{-1-b},\ b>0$).  $\Omega$ is a bounded domain
satisfying (H).  Then the function $V(x,t)$ is an upper solution to
problem (\ref{e:critical}) and exists for all $t>0$.
\end{thm}

In order to prove Theorem~\ref{theo:critical}, we first derive a
number of preliminary facts on $d(x)$.

\begin{lem}\label{lem:critical1}
Assume $x_0=(x_{10},x_{20},\cdots,x_{n0})$, $\Omega_i=\{x\in
R^n:|x-x_0|<R_i\}, i=1,2$ and $R_1>R_2$. Let
$d(x)=dist(x,\partial\Omega_1), x\in\Omega_1\setminus\Omega_2$. Then
$\Delta d(x)=(1-n)/(|x-x_0|)$.
\end{lem}

\begin{lem}
$\Omega$ is a bounded domain satisfying (H).  Then there exists
$\varepsilon>0$ such that $\Delta d\leq0$ for $x\in
\Omega_\varepsilon$.
\end{lem}

\begin{proof}
Here we only consider the case of $n=2$.  As for $n=1$ or $n\geq3$,
the proof is completely similar.  Divided $\partial\Omega$ into $m$
parts and taking $m$ large enough such that the largest arc is
sufficiently small.  Let $A_1,A_2,\cdots,A_m$ be the division
points. For any arc $\widehat{{A_iA}}_{i+1}(1\leq i\leq m-1)$,
choosing $C\in\widehat{{A_iA}}_{i+1}$ such that
$|\widehat{{A_iC}}|=|\widehat{{CA}}_{i+1}|$.  By the definition of
$\Omega$, there exists a circle $\Omega_1=\{x\in R^2:|x-x_0|<R_1\}$
such that $\Omega_1(\Omega\subset\Omega_1)$ is tangent to $\Omega$
at the point $C$.Taking
$A_i^\prime,A_{i+1}^\prime\in\partial\Omega_1$ such that the
segments $A_i^\prime x_0,A_{i+1}^\prime x_0$ intersect
$\partial\Omega$ at $A_i,A_{i+1}$, respectively. Since
$\widehat{{A_iA}}_{i+1}$ is sufficiently small, we have
$\widehat{A_iA}_{i+1}\sim \widehat{{A_i^{'}A^{'}}}_{i+1}$. From
Lemma~\ref{lem:critical1}, there exists a constant
$\varepsilon_{\widehat{{A_iA}}_{i+1}}>0$ such that
\[
  \Delta d(x)\leq\frac{-1}{2|x-x_0|}<0,\ \ \
  x\in\Big\{x\in\Omega:d(x,\widehat{{A_iA}}_{i+1})<\varepsilon_{\widehat{{A_iA}}_{i+1}}\Big\}.
\]
Set
$\varepsilon=\min\{\varepsilon_{\widehat{{A_iA}}_{i+1}},\varepsilon_{\widehat{{A_1A}}_{m}},i=1,2,\cdots,m-1\}$.
Then
\[
  \Delta d(x)\leq0, \ \ \
  x\in\Omega_\varepsilon=\{x\in\Omega:d(x,\partial\Omega)<\varepsilon\}.
\]
The proof is completed.
\end{proof}

Now we give the proof of Theorem~\ref{theo:critical}.

\begin{proof}
Case1: We assume $f(s)$ to be that
$\liminf_{s\rightarrow\infty}g(s)>C>0$ and
$\lim_{s\rightarrow\infty}\mu(s)f(s)=C_0>0$.  Then taking into
account the relation (\ref{e:upper12}), for $d(x)\geq
\varepsilon(t)$, we get
\begin{equation*}
\begin{split}
\mathcal {F}(V)&=V_t-\Delta
V-\frac{2|\partial\Omega|^2f(V)}{(\int_\Omega f(V)dx)^2}\\
&\geq
\dot{M}(t)-\frac{2|\partial\Omega|^2f(M)}{(|\Omega\setminus\Omega_\varepsilon|f(M)+\varepsilon|\partial\Omega|\sqrt{\frac{2}{\mu}}\sqrt{\int_0^M
  f(s)ds}-\varepsilon^2|\partial\Omega|
  K\frac{M}{\mu})^2}\\
  &\geq\dot{M}(t)-\frac{\mu(M)f(M)}{\varepsilon^2(\frac{|\Omega\setminus\Omega_\varepsilon|\sqrt{\mu}f(M)}{\sqrt{2}\varepsilon|\partial\Omega|}+\int_0^M
  f(s)ds-\frac{K\varepsilon M}{\sqrt{2\mu}})^2}\\
  &\geq \dot{M}(t)-\frac{\mu(M)f(M)}{\varepsilon^2(\frac{|\Omega|\sqrt{\mu}f(M)}{2\sqrt{2}\varepsilon|\partial\Omega|}+\int_0^M
  f(s)ds-\frac{K\varepsilon M}{\sqrt{2\mu}})^2}\ \ \ \ \ \ \  {\rm for}\ \
  \varepsilon(t)\ll1.
\end{split}
\end{equation*}
Choosing $K_1=(C_0|\Omega|)/(8K|\partial\Omega|)$ and
$\varepsilon(t)=(K_1/M)^{1/2}$, we have $0<\varepsilon(M)\ll1$ for
$M\gg1$.  Moreover, from (\ref{e:upper10}), (\ref{e:upper11}) and
$\lim_{M\rightarrow\infty}\mu(M)f(M)=C_0$, we obtain
\begin{equation*}
\begin{split}
&\frac{|\Omega|\sqrt{\mu}f(M)}{2\sqrt{2}\varepsilon|\partial\Omega|}+\int_0^M
  f(s)ds-\frac{K\varepsilon M}{\sqrt{2\mu}}\\
  &\geq \frac{|\Omega|\sqrt{\mu}f(M)}{2\sqrt{2}\varepsilon|\partial\Omega|}+\int_0^M
  f(s)ds-\frac{\sqrt{2}K_1K\sqrt{\mu}f(M)}{C_0\varepsilon}\\
&=\frac{|\Omega|\sqrt{\mu}f(M)}{4\sqrt{2}\varepsilon|\partial\Omega|}+\int_0^M
  f(s)ds\ \ \ {\rm for}\ \ M\gg1.
\end{split}
\end{equation*}
Since
\[
  \frac{|\Omega|\sqrt{\mu}f(M)}{4\sqrt{2}\varepsilon|\partial\Omega|
  F(M)}=\frac{|\Omega|\sqrt{\mu}f(M)}{4\sqrt{2}\varepsilon|\partial\Omega| (1-\int_0^M
  f(s)ds)}\geq\frac{|\Omega|C}{4\sqrt{2}\varepsilon|\partial\Omega|}>1 \ \ \ {\rm for}\ \
  M\gg1,
\]
which implies
\begin{equation}\label{e:upper13}
  \frac{|\Omega|\sqrt{\mu}f(M)}{4\sqrt{2}\varepsilon|\partial\Omega|}+\int_0^M
  f(s)ds>1\ \ \ {\rm for}\ \ M\gg1.
\end{equation}
Taking $M(t)$ to satisfy
\begin{equation}\label{e:upper14}
  \dot{M}(t)=\frac{\mu(M)f(M)}{\varepsilon^2(M)},\ \ \ t>0,
\end{equation}
we obtain
\[
  \mathcal {F}(V)>\dot{M}(t)-\frac{\mu(M)f(M)}{\varepsilon^2(M)}=0\
  \ \ {\rm for} \ d(x)\geq\varepsilon(t)(x\in\Omega)\ \ {\rm and}\ \
  M\gg1.
\]
By integrating (\ref{e:upper14}), we have
\[
  \int_{M(0)}^{M(t)}\frac{\varepsilon^2(s)}{\mu(s)f(s)}ds=t,
\]
and taking into account $\lim_{s\rightarrow\infty}\mu(s)f(s)=C_0$,
we obtain
\[
  \frac{K_1}{1+C_0}\int_{M(0)}^{M(t)}\frac{1}{s}ds<t\ \ \ {\rm for}\ \
  M(0)\gg1.
\]
The last inequality implies that if $M(t)\rightarrow\infty$ then
$t\rightarrow\infty$.

Also for $0<d(x)\leq\varepsilon(t)(x\in\Omega)$, we have
\begin{equation*}
\begin{split}
\mathcal
{F}(V)&=w_{\mu}\big(y(x,t);\mu(t)\big)\dot{\mu}(t)+\frac{dw\big(y(x,t);\mu(t)\big)}{dy}\dot{y}(t)-\Delta
w-\frac{2|\partial\Omega|^2f(w)}{(\int_\Omega f(V)dx)^2}\\
&=w_{\mu}\big(y(x,t);\mu(t)\big)\dot{\mu}(t)-\frac{dw\big(y(x,t);\mu(t)\big)}{dy}\frac{d(x)}{\varepsilon^2}\dot{\varepsilon}(t)-\frac{\Delta
d}{\varepsilon}\frac{dw}{dy}+\frac{\mu}{\varepsilon^2}f(w)-\frac{2|\partial\Omega|^2f(w)}{(\int_\Omega
f(V)dx)^2}.
\end{split}
\end{equation*}
Since $w_{\mu}>0,\ \dot{\mu}(t)>0,\ \dot{\varepsilon}(t)<0$,
$dw/dy\geq0$ and $\Delta d(x)\leq0$ for $M\gg1$, we have
\begin{equation*}
\mathcal {F}(V)\geq\frac{\mu
f(w)}{\varepsilon^2}-\frac{2|\partial\Omega|^2f(w)}{(\int_\Omega
f(V)dx)^2}\geq\frac{\mu
f(w)}{\varepsilon^2}\Big(1-\frac{1}{(\frac{|\Omega|\sqrt{\mu}f(M)}{4\sqrt{2}\varepsilon|\partial\Omega|}+\int_0^M
  f(s)ds)^2}\Big)>0\ \ \ {\rm for}\ \
  M\gg1.
\end{equation*}

Case2: Now let $f$ be such
$\liminf_{s\rightarrow\infty}\mu(s)f(s)/s=C_1>0(C_1\leq2)$ and
$\liminf_{s\rightarrow\infty}g(s)>C>0$.  For this case it is enough
to consider $\varepsilon(t)$ to be constant such that $\Delta d\leq
0$ for $x\in\Omega_\varepsilon$. Moreover,  we choose $\varepsilon$
to satisfy
\[
 \frac{|\Omega\setminus\Omega_\varepsilon|}{\sqrt{2}\varepsilon|\partial\Omega|}-\frac{\sqrt{2}K\varepsilon}{C_1}>\frac{1}{C}.
\]
For $d(x)\geq\varepsilon(x\in\Omega)$, we have
\begin{equation*}
\begin{split}
\mathcal {F}(V) &\geq
\dot{M}(t)-\frac{2|\partial\Omega|^2f(M)}{(|\Omega\setminus\Omega_\varepsilon|f(M)+\varepsilon|\partial\Omega|\sqrt{\frac{2}{\mu}}\sqrt{\int_0^M
  f(s)ds}-\varepsilon^2|\partial\Omega|
  K\frac{M}{\mu})^2}\\
  &\geq\dot{M}(t)-\frac{\mu(M)f(M)}{\varepsilon^2(\frac{|\Omega\setminus\Omega_\varepsilon|\sqrt{\mu}f(M)}{\sqrt{2}\varepsilon|\partial\Omega|}+\int_0^M
  f(s)ds-\frac{K\varepsilon M}{\sqrt{2\mu}})^2}\\
  &\geq \dot{M}(t)-\frac{\mu(M)f(M)}{\varepsilon^2(\frac{|\Omega\setminus\Omega_\varepsilon|\sqrt{\mu}f(M)}{\sqrt{2}\varepsilon|\partial\Omega|}+\int_0^M
  f(s)ds-\frac{\sqrt{2}K\varepsilon f(M)\sqrt{\mu}}{C_1})^2}\\
  &\geq \dot{M}(t)-\frac{f(M)\mu(M)}{\varepsilon^2(\frac{f(M)\sqrt{\mu(M)}}{C}+\int_0^M
  f(s)ds)^2}\ \ \  {\rm for}\ \ \ M\gg1.
\end{split}
\end{equation*}
Since
\[
  \frac{f(M)\sqrt{\mu(M)}}{C F(M)}=\frac{f(M)\sqrt{\mu(M)}}{C(1-\int_0^M
  f(s)ds)}>\frac{1}{C}C=1,
\]
which implies
\[
  \frac{f(M)\sqrt{\mu(M)}}{C }+\int_0^M
  f(s)ds>1.
\]
Hence $\mathcal {F}(V)>0$ for $d(x)\geq\varepsilon$ and $M\gg1$,
provided that $M(t)$ satisfies
\begin{equation}\label{e:upper16}
  \dot{M}(t)=\frac{\mu(M)f(M)}{\varepsilon^2},\ \ \ t>0.
\end{equation}
By integrating (\ref{e:upper16}), we have
\[
  \int_{M(0)}^{M(t)}\frac{\varepsilon^2}{\mu(s)f(s)}ds=t,
\]
and taking into account (\ref{e:upper9}) we obtain
\[
  \frac{\varepsilon^2}{2}\int_{M(0)}^{M(t)}\frac{1}{s}ds\leq t,
\]
which implies that if $M(t)\rightarrow\infty$ then
$t\rightarrow\infty$.

For $0\leq d(x)\leq\varepsilon$, we have
\begin{equation*}
\begin{split}
\mathcal {F}(V)
&=w_{\mu}\big(y(x,t);\mu(t)\big)\dot{\mu}(t)-\frac{\Delta
d}{\varepsilon}\frac{dw}{dy}+\frac{\mu}{\varepsilon^2}f(w)-\frac{2|\partial\Omega|^2f(w)}{(\int_\Omega
f(V)dx)^2}\\
 &\geq \frac{\mu(M) f(w)}{\varepsilon^2}-\frac{f(w)\mu(M)}{\varepsilon^2(\frac{f(M)\sqrt{\mu(M)}}{C }+\int_0^M
  f(s)ds)^2}>0\ \ \  {\rm for}\ \ \ M\gg1.
\end{split}
\end{equation*}
Therefore, we finally get that, in each case, $V(x,t)$ is an upper
solution to problem (\ref{e:critical}) for all $t>0$.  The proof is
completed.
\end{proof}

Thus we formulate this main result of this Section in the following
theorem.

\begin{thm}
If $f(s)$ satisfies the hypotheses of Theorem~\ref{theo:critical},
and $\Omega$ satisfies (H), then $u(x,t)$ is a global-in-time
solution to problem (\ref{e:critical}) and $u(x,t)\rightarrow\infty$
as $t\rightarrow\infty$, i.e. $u(x,t)$ diverges globally in
$\Omega$.
\end{thm}

\section{Asymptotic behavior of the blow-up solutions}
In this Section, we deal with the blow-up solutions of problem
(\ref{e:main}).

\begin{thm}\label{theo:blowup}
Let $f(s)$ satisfy (\ref{e:f}), $\int_0^\infty f(s)ds=1$, $p=2$ and
$\Omega$ satisfy (H).  If $\lambda>\lambda^*=2|\partial\Omega|^2$,
the solution of the problem (\ref{e:main}) blows up globally in
finite time $T$ .
\end{thm}
\begin{proof}
By Theorem~\ref{theo:steady1}, in the case of
$\lambda>\lambda^*=2|\partial\Omega|^2$ and $\int_0^\infty
f(s)ds=1$, there is no steady solution to (\ref{e:steady2}).  Since
$\lambda(\mu)<\lambda$ for any $\mu>0$, we can find an increasing
lower solution $v(x,t)=w(x;\mu(t))$ with $\mu$ and
$v\rightarrow\infty$ as $t\rightarrow T\leq\infty$.  Thus $u(x,t)$
is globally unbounded.  We shall show that $T<\infty$.  Therefore,
we look for a lower solution $V(x,t)$ which blows up at a finite
time ( $V(x,t)$ satisfy (\ref{e:upper})$-$(\ref{e:upper4}) ).  From
(\ref{e:upper4}) and (\ref{e:upper6}), we have
\begin{equation*}
\begin{split}
  \int_{\Omega}f(V)dx&=\int_{\Omega\setminus\Omega_\varepsilon}f(M)dx+\int_{\Omega_\varepsilon}f(w)dx\leq|\Omega|f(M)
  +|\partial\Omega|\int_{\delta(t)}^1f\big(w(x_1,0');\mu(t)\big)dx_1\\
  &\leq|\Omega|f(M)+|\partial\Omega|\varepsilon\sqrt{\frac{2}{\mu}}
  =\sqrt{2}|\partial\Omega|f(M)(\frac{|\Omega|}{\sqrt{2}|\partial\Omega|}+\alpha),
\end{split}
\end{equation*}
on choosing $\alpha=\varepsilon/(\sqrt{\mu}f(M))$, where $\alpha$ is
a suitable chosen constant; in particular choose
$\alpha>|\Omega|/(\sqrt{\lambda}-\sqrt{2}|\partial\Omega|)$ for
$\lambda>\lambda^*=2|\partial\Omega|^2$.  Such an $\alpha$ gives
\[
  3\Lambda=\frac{\lambda}{(|\Omega|+\sqrt{2}\alpha|\partial\Omega|)^2}-\frac{1}{\alpha^2}>0.
\]
From (\ref{e:upper10}), we also note that with such a fixed
$\alpha$, $\varepsilon\rightarrow0$ as $M\rightarrow\infty$.
Integrating (\ref{e:upper7}) on $(0,r)$, we get
\begin{equation}\label{e:lower}
  \int_0^w\frac{ds}{\sqrt{F(s)-F(M)}}=\frac{\sqrt{2\mu}r}{\varepsilon}=\frac{\sqrt{2}r}{\alpha
  f(M)}.
\end{equation}

For $x\in\Omega\setminus\Omega_\varepsilon$,
\begin{equation*}
\mathcal {F}(V) =\dot{M}-\frac{\lambda f(M)}{(
\int_{\Omega}f(V)dx)^2}\leq\dot{M}-\frac{\lambda}{2|\partial\Omega|^2f(M)(\frac{|\Omega|}{\sqrt{2}|\partial\Omega|}+\alpha)^2}
 \leq \dot{M}-\frac{\Lambda}{f(M)}\leq0,
\end{equation*}
on choosing $\dot{M}\leq\Lambda/f(M)$.

For $x\in\Omega_\varepsilon$, we first differentiate (\ref{e:lower})
with respect to $t$ and get
\begin{equation*}
\begin{split}
  w_t&=-\frac{f'(M)}{f(M)}\dot{M}(t)[F(w)-F(M)]^{\frac{1}{2}}\int_0^w\frac{ds}{\sqrt{F(s)-F(M)}}\\
  &+\frac{1}{2}f(M)\dot{M}(t)[F(w)-F(M)]^{\frac{1}{2}}\int_0^w[F(s)-F(M)]^{-\frac{3}{2}}ds:=A+B.
\end{split}
\end{equation*}
For $A$, from (\ref{e:upper00}) we have
\begin{equation*}
\begin{split}
 A&=-\frac{f'(M)}{f(M)}\dot{M}(t)[F(w)-F(M)]^{\frac{1}{2}}\int_0^w\frac{ds}{\sqrt{F(s)-F(M)}}\\
  &\leq-\frac{2f'(M)}{f^{\frac{3}{2}}(M)}M\dot{M}(t)f^{\frac{1}{2}}(w)\leq\frac{\Lambda
  f(w)}{f^2(M)},
\end{split}
\end{equation*}
provided that
\[
  \dot{M}(t)\leq-\frac{\Lambda}{2M f'(M)}
\]
and taking into account that $f'(s)\leq0$ so that $f(w)/f(M)\geq1$
for $w\leq M$.  For $B$ we have
\begin{equation*}
\begin{split}
 B&=\frac{1}{2}f(M)\dot{M}(t)[F(w)-F(M)]^{\frac{1}{2}}\int_0^w[F(s)-F(M)]^{-\frac{3}{2}}ds\\
  &\leq\frac{f^{\frac{1}{2}}(w)}{f^{\frac{1}{2}}(M)}\dot{M}(t)\leq\frac{\Lambda
  f(w)}{f^2(M)},
\end{split}
\end{equation*}
provided that
\[
  \dot{M}(t)\leq\frac{\Lambda}{ f(M)}.
\]
Also, using (\ref{e:upper3}) and (\ref{e:upper7}), we have the
estimate
\begin{equation*}
\begin{split}
 -\Delta w&=-w_r\Delta d +\frac{\mu}{\varepsilon^2}f(w)\leq K w_r+\frac{\mu}{\varepsilon^2}f(w)\ \ ({\rm Using} \ |\Delta d|\leq K)\\
  &=\frac{K\sqrt{2\mu}}{\varepsilon}[F(w)-F(M)]^{\frac{1}{2}}+\frac{f(w)}{\alpha^2f^2(M)}\\
  &\leq\frac{\sqrt{2}K}{\alpha}\frac{(Mf(M))^{\frac{1}{2}}f(w)}{f^2(M)}+\frac{f(w)}{\alpha^2f^2(M)}\\
  &\leq \frac{\Lambda f(w)}{f^2(M)}+\frac{f(w)}{\alpha^2f^2(M)},\ \
  \  {\rm for}\ M\gg1,
\end{split}
\end{equation*}
since $Mf(M)\rightarrow0$ as $M\rightarrow\infty$.  Thus for
$x\in\Omega_\varepsilon$ if
\begin{equation}\label{e:lower1}
  0\leq\dot{M}(t)=\min\{\frac{\Lambda}{ f(M)},\ -\frac{\Lambda}{2M
  f'(M)}\}
\end{equation}
and using the previous estimate we obtain
\begin{equation*}
\begin{split}
  \mathcal {F}(V)&=w_t-\Delta w-\frac{\lambda f(w)}{(\int_\Omega
  f(V)dx)^2}\\
  &=A+B-w_r\Delta d +\frac{\mu}{\varepsilon^2}f(w)-\frac{\lambda f(w)}{(\int_\Omega f(V)dx)^2}\\
  &\leq\frac{3\Lambda
  f(w)}{f^2(M)}+\frac{f(w)}{\alpha^2f^2(M)}-\frac{\lambda
  f(w)}{2|\partial\Omega|^2f^2(M)(\frac{|\Omega|}{\sqrt{2}|\partial\Omega|}+\alpha)^2}=0.
\end{split}
\end{equation*}
Also $V(x,t)=u(x,t)=0$ on the boundary $\partial\Omega$ and taking
$V(x,0)\leq u_0(x)$, the function $V(x,t)$ is a lower solution to
the problem (\ref{e:main}). Hence $u(x,t)\geq V(x,t)$ for $M$ is
large enough (after some time at which $u(x,t)$ is sufficiently
large if $T=\infty$).

Now we show that $u(x,t)$ blows up in finite time.  Indeed, from
(\ref{e:lower1}) we have
\[
  \Lambda\frac{dt}{dM}=\max\{ f(M),\ -2M f'(M)\}\leq f(M)-2M f'(M) \
  \ \ (f'(s)\leq0)
\]
or
\[
  \Lambda t\leq \int^M[f(s)-2sf'(s)]ds<\infty,
\]
since $Mf(M)\rightarrow0$ as $M\rightarrow\infty$ and $\int_0^\infty
f(s)ds=1$.  Hence $V(x,t)$ blows up at $t^*<\infty$ and $u(x,t)$
must blow up at $T\leq t^*<\infty$.

As for the blow-up is global from the fact
\[
  \int_\Omega f(u)dx\rightarrow0\ \ \  {\rm as}\ \ t\rightarrow T.
\]
Indeed,
\[
  \dot{M}\leq\frac{\lambda f(M)}{(\int_\Omega f(u)dx)^2}=h(t),
\]
giving
\[
  M(t)-M(0)\leq\int_0^t h(s)ds\rightarrow\infty\ \ \  {\rm as}\ \ t\rightarrow
  T.
\]
This implies $\int_\Omega f(u)dx\rightarrow0$ as $ t\rightarrow T$
since $f(s)$ is bounded.  Thus, for
$\lambda>\lambda^*=2|\partial\Omega|^2$, $u(x,t)$ blows up globally.
The proof is completed.
\end{proof}

\begin{thm}
Let $f(s)$ satisfy (\ref{e:f}), $\int_0^\infty f(s)ds=1$, $p>2$ and
$\Omega$ satisfy (H).  Then there exists a critical value
$\lambda^*$ such that for $\lambda>\lambda^*$ or for any
$0<\lambda\leq\lambda^*$ but with initial data sufficiently large,
the solution of the problem (\ref{e:main}) blows up globally in
finite time $T$.
\end{thm}
\begin{proof}
Using Theorem~\ref{theo:steady1}, we know that for
$\lambda>\lambda^*$ or for any $0<\lambda\leq\lambda^*$ but with
initial data $u_0$ more than the greater steady state $u(x,t)$ is
globally unbounded(see\cite{L2}). In order to prove $u(x,t)$ blows
up in finite time $T<\infty$, we also look for a lower solution
$V(x,t)$ to satisfy (\ref{e:upper})$-$(\ref{e:upper4}).  Then
\begin{equation*}
\begin{split}
  \int_{\Omega}f(V)dx&=\int_{\Omega\setminus\Omega_\varepsilon}f(M)dx+\int_{\Omega_\varepsilon}f(w)dx\leq|\Omega|f(M)
  +|\partial\Omega|\int_{\delta(t)}^1f\big(w(x_1,0);\mu(t)\big)dx_1\\
  &\leq|\Omega|f(M)+|\partial\Omega|\varepsilon\sqrt{\frac{2}{\mu}}
  =\sqrt{2}|\partial\Omega|f(M)(\frac{|\Omega|}{\sqrt{2}|\partial\Omega|}+1),
\end{split}
\end{equation*}
on choosing $\varepsilon=\sqrt{\mu}f(M)$.  From (\ref{e:upper10}),
we also note that $\varepsilon\rightarrow0$ as $M\rightarrow\infty$.

For $x\in\Omega\setminus\Omega_\varepsilon$,
\begin{equation*}
\begin{split}
\mathcal {F}(V) &=\dot{M}-\frac{\lambda f(M)}{(
\int_{\Omega}f(V)dx)^p}\leq\dot{M}-\frac{\lambda}{(\sqrt{2}|\partial\Omega|)^p
f^{p-1}(M)(\frac{|\Omega|}{\sqrt{2}|\partial\Omega|}+1)^p} \\
& \leq \dot{M}-\frac{1}{f(M)}\leq0 \ \ \  {\rm for}\  M\gg1,
\end{split}
\end{equation*}
on choosing $\dot{M}\leq1/f(M)$ and taking into account $p>2$ and
$f(M)\rightarrow0$ as $M\rightarrow\infty$.

For $x\in\Omega_\varepsilon$, similar to the proof of
Theorem~\ref{theo:blowup}, we have $w_t=A+B$.  For $A$, from
(\ref{e:upper00}) we have
\begin{equation*}
\begin{split}
 A&=-\frac{f'(M)}{f(M)}\dot{M}(t)[F(w)-F(M)]^{\frac{1}{2}}\int_0^w\frac{ds}{\sqrt{F(s)-F(M)}}\\
  &\leq-\frac{2f'(M)}{f^{\frac{3}{2}}(M)}M\dot{M}(t)f^{\frac{1}{2}}(w)\leq\frac{
  f(w)}{f^2(M)},
\end{split}
\end{equation*}
provided that
\[
  \dot{M}(t)\leq-\frac{1}{2M f'(M)}.
\]
For $B$ we have
\begin{equation*}
\begin{split}
 B&=\frac{1}{2}f(M)\dot{M}(t)[F(w)-F(M)]^{\frac{1}{2}}\int_0^w[F(s)-F(M)]^{-\frac{3}{2}}ds\\
  &\leq\frac{f^{\frac{1}{2}}(w)}{f^{\frac{1}{2}}(M)}\dot{M}(t)\leq\frac{
  f(w)}{f^2(M)},
\end{split}
\end{equation*}
provided that
\[
  \dot{M}(t)\leq\frac{1}{ f(M)}.
\]
Also, using (\ref{e:upper3}) and (\ref{e:upper7}), we have the
estimate
\begin{equation*}
\begin{split}
 -\Delta w&=-w_r\Delta d +\frac{\mu}{\varepsilon^2}f(w)\leq K w_r+\frac{\mu}{\varepsilon^2}f(w)\ \ ({\rm Using} \ |\Delta d|\leq K)\\
  &=\frac{K\sqrt{2\mu}}{\varepsilon}[F(w)-F(M)]^{\frac{1}{2}}+\frac{f(w)}{f^2(M)}\\
  &\leq\sqrt{2}K\frac{(Mf(M))^{\frac{1}{2}}f(w)}{f^2(M)}+\frac{f(w)}{f^2(M)}\\
  &\leq \frac{2f(w)}{f^2(M)},\ \
  \  {\rm for}\ M\gg1,
\end{split}
\end{equation*}
since $Mf(M)\rightarrow0$ as $M\rightarrow\infty$.  Thus for
$x\in\Omega_\varepsilon$ if
\begin{equation*}
  0\leq\dot{M}(t)=\min\{\frac{1}{ f(M)},\ -\frac{1}{2M
  f'(M)}\}
\end{equation*}
and using the previous estimate we obtain
\begin{equation*}
\begin{split}
  \mathcal {F}(V)&=w_t-\Delta w-\frac{\lambda f(w)}{(\int_\Omega
  f(V)dx)^p}\\
  &=A+B-w_r\Delta d+\frac{\mu}{\varepsilon^2}f(w)-\frac{\lambda f(w)}{(\int_\Omega f(V)dx)^p}\\
  &\leq\frac{4
  f(w)}{f^2(M)}-\frac{\lambda f(w)}{(2|\partial\Omega|)^p f^p(M)(\frac{|\Omega|}{\sqrt{2}|\partial\Omega|}+1)^p}\leq0\ \
  \  {\rm for}\ M\gg1,
\end{split}
\end{equation*}
since $p>2$ and $f(M)\rightarrow0$ as $M\rightarrow\infty$.

Also $V(x,t)=u(x,t)=0$ on the boundary $\partial\Omega$ and taking
$V(x,0)\leq u_0(x)$, the function $V(x,t)$ is a lower solution to
the problem (\ref{e:main}). Hence $u(x,t)\geq V(x,t)$ for $M$ is
large enough (after some time at which $u$ is sufficiently large if
$T=\infty$).

The rest proof is same as the Theorem~\ref{theo:blowup}, so we omit
it here.
\end{proof}

Now we will consider the Dirichlet problem, which we rewrite
(\ref{e:main}) as
\begin{eqnarray*}
\begin{cases}
  u_{t}=\Delta u+g(t)f(u),\ \ \ &x\in \Omega,\ t>0,\\
  u(x,t)=0,            &x\in\partial \Omega,\ t>0,\\
  u(x,0)=u_{0}(x),     &x\in \Omega,
\end{cases}
\end{eqnarray*}
where $g(t)=\lambda/(\int_\Omega f(u)dx)^p$, $\Omega$ is defined as
Theorem~\ref{theo:steady1}.

We seek a formal asymptotic approximation for $u(x,t)$ near the
blow-up time $T$, still taking $f$ to be decreasing and to satisfy
$\int_0^\infty f(s)ds=1$.  Set $M(t)=\max_{x\in\Omega}u(x,t)$.

As in \cite{KT}, we obtain that $\lim_{t\rightarrow T}g(t)=\infty$
and $u(x,t)\sim M$ except in some boundary layers near
$\partial\Omega$. In the main core(outer) region we neglect $\Delta
u$, so
\[
  \frac{dM}{dt}\sim g(t)f(M)
\]
and significant contributions to the integral $\int_\Omega f(u)dx$
can come from the largest(core) region which has volume $\sim
|\Omega|$(contribution$\sim|\Omega|f(M)$) and from the boundary
layers where $f$ is large, $f(u)$ is $O(1)$ where $u(x,t)$ is
$O(1)$.  If the boundary layers have volume $O(\delta)$, for some
small $\delta$, then to obtain a balance involving $\Delta u$,
either $\delta^{-2}=O(g)$ or $\delta^{-2}=O\big((T-t)^{-1}\big)$,
whichever is the larger, see \cite{L2}.

Supposing that $g(t)\ll(T-t)^{-1}$ for $t\rightarrow T$ the
contribution to the the integral from the boundary layer is
$O(\delta)=O\big((T-t)^{\frac{1}{2}}\big)$, whereas
\[
  \int_\Omega f(u)dx=O\big(g(t)^{-\frac{1}{p}}\big)\gg
  (T-t)^{\frac{1}{p}}\geq (T-t)^{\frac{1}{2}}\ \ \  {\rm as}\ \ \ t\rightarrow T.
\]
This suggests that the core dominates and
\[
  \int_\Omega f(u)dx\sim |\Omega| f(M).
\]
Then
\[
  g(t)\sim \frac{\lambda}{|\Omega|^p f^p(M)},\ \ \  f(M)\sim
  \frac{1}{|\Omega|}(\frac{\lambda}{g})^\frac{1}{p},
\]
and
\[
  \frac{dM}{dt}\sim
  g(t)f(M)\sim\frac{1}{|\Omega|}\lambda^{\frac{1}{p}}g^{\frac{p-1}{p}}\ll(T-t)^{\frac{1-p}{p}}\
  \ \  {\rm for}\ \ \ t\rightarrow T.
\]
This would indicate that $M$ is actually bounded as $t\rightarrow
T$.  Contradicting the occurrence of blow-up.

Next we suppose that $g(t)=O\big((T-t)^{-1}\big)$ for $t\rightarrow
T$.  Since
\[
  |\Omega|f(M)\lesssim\int_\Omega
  f(u)dx=(\frac{\lambda}{g})^\frac{1}{p},
\]
we must have $f(M)\leq O\big((T-t)^{\frac{1}{p}}\big)$.  Again,
\[
  \frac{dM}{dt}\sim
  g(t)f(M)\leq O\big((T-t)^{\frac{1-p}{p}}\big),
\]
which contradicts the assumption of blow-up.  There remains only one
possibility:
\[
  g(t)\gg(T-t)^{-1}\ \ \  {\rm for}\ \ \ t\rightarrow T.
\]

The boundary layer has volume
$O\big(g(t)^{-\frac{1}{2}}\big)\ll(T-t)^{\frac{1}{2}}$, where
$u(x,t)$ is $O(1)$ and $u_t$ is negligible compare to $\Delta u$.
There has to be a balance between $\Delta u$ and $g(t)f(u)$, that is
,
\[
  -\Delta u\sim g(t)f(u).
\]
Without loss of generality, we assume that the hyperplane $\{x\in
R^n:x_1=1\}$ is tangent to $\Omega$ at $y_0(y_0=(1,0'))$, and
$\Omega$ lies in the half-space $\{x:x_1<1\}$. Writing
$x_1=1-y/\sqrt{g}(y/\sqrt{g}\ll1)$ gives
\begin{eqnarray}\label{e:general}
\begin{cases}
  -u_{yy}\big((y,0')\big)\sim f\big(u((y,0'))\big),\ \ \ &y>0,\\
  u\big((y,0')\big)=0,\ \ \ &y=0,\\
  u\big((y,0')\big)\gg1\gg u_y\big((y,0')\big),\ \ \ &y\gg1.
\end{cases}
\end{eqnarray}
Multiplying both sides of (\ref{e:general}) by $u_y\big((y,0')\big)$
and integrating, we get
\[
  u_y^2\big((y,0')\big)\sim 2F\big(u((y,0'))\big),
\]
where $F\big(u((y,0'))\big)=\int_{u((y,0'))}^{\infty}f(s)ds$.
Integrating again gives $u\big((y,0')\big)\sim U(y)$, where
\begin{equation}\label{e:blow}
  \sqrt{2}y=\int_0^{U(y)}F^{-\frac{1}{2}}(s)ds.
\end{equation}
Since $y_0$ is arbitrary, it follows from (\ref{e:blow}) that the
boundary layers contribute to a total amount
\[
 \int_{d(x,\partial\Omega)\leq y/\sqrt{g}}f(u)dx\sim
 |\partial\Omega|\int_{x_1}^1f\big(u((x_1,0'))\big)dx_1\sim\frac{|\partial\Omega|}{\sqrt{g}}\int_0^\infty
 f\big(U(y)\big)dy,
\]
this is automatically of the correct size $g(t)=\lambda/(\int_\Omega
f(u)dx)^p$.  It should also be observed that
\[
  \int_0^\infty f\big(U(y)\big)dy=U'(0).
\]

Looking at the following steady problem
\[
  w''+\mu f(w)=0,\ \ \ -1<x<1;\ \ \ w(\pm1)=0.
\]
Set $M(\mu)=\max_{-1<x<1}w(x)=w(0)$ and $x=1-y/\sqrt{\mu}$, then
\[
  \frac{d^2w}{dy^2}+f(w)=0, \ \ \ w(0)=0,\ \ \
  \frac{dw}{dy}|_{y=\sqrt{\mu}}=0,\ \ \  w(\sqrt{\mu})=M.
\]
From Lemma~\ref{prop:4} (in case of $n=1$), we have
\[
  \lim_{\mu\rightarrow\infty}-\frac{1}{\sqrt{\mu}}\frac{dw(1)}{dx}=\sqrt{2}.
\]
which implies
\[
  \lim_{\mu\rightarrow\infty}\frac{dw(0)}{dy}=\sqrt{2},
\]
and it appears that the problem in limit of large $\mu$ is the same
as the asymptotic problem (\ref{e:general}).  Thus,
\[
  \int_0^\infty
  f\big(U(y)\big)dy=U'(0)=\lim_{\mu\rightarrow\infty}\frac{dw(0)}{dy}=\sqrt{2}.
\]
We deduce that the contribution to $\int_\Omega f(u)dx$ from the
boundary layers$\sim \sqrt{2}|\partial\Omega|/\sqrt{g}$.

Now
\[
  \int_\Omega f(u)dx\sim|\Omega|f(M)+\sqrt{2}|\partial\Omega|/\sqrt{g}
\]
and
\[
  g\sim\frac{\lambda}{(|\Omega|f(M)+\sqrt{2}|\partial\Omega|/\sqrt{g})^p}\
  \ \ {\rm for}\ \  t\rightarrow T (g,M\rightarrow\infty).
\]
We see that
\[
  \lambda^{\frac{1}{p}}\sim
  g^{\frac{1}{p}}(|\Omega|f(M)+\sqrt{2}|\partial\Omega|/\sqrt{g})=|\Omega|f(M)g^{\frac{1}{p}}+\sqrt{2}|\partial\Omega|
  g^{\frac{2-p}{2p}},
\]
i.e.,
\begin{enumerate}
    \item[(i)] If $p=2$, then
    $f(M)\sim\displaystyle\frac{\sqrt{\lambda}-\sqrt{2}|\partial\Omega|}{|\Omega|\sqrt{g}}$.
    \item[(ii)] If $p>2$, then $f(M)\sim
    \displaystyle\frac{1}{|\Omega|}(\frac{\lambda}{g})^{\frac{1}{p}}$.
\end{enumerate}
Therefore, in the core region $u(x,t)\sim M$ which satisfies
\begin{equation}\label{e:blow1}
  \frac{dM}{dt}\sim
  g(t)f(M)\sim\frac{\Lambda_1^2}{f(M)}\  \ \ {\rm if} \ \  p=2,
\end{equation}
where
$\Lambda_1=(\sqrt{\lambda}-\sqrt{2}|\partial\Omega|)/|\Omega|$, and
\begin{equation}\label{e:blow2}
  \frac{dM}{dt}\sim
  g(t)f(M)\sim\frac{\Lambda_2}{f^{p-1}(M)}\  \ \ {\rm if} \ \  p>2,
\end{equation}
where $\Lambda_2=\lambda/|\Omega|^p$.

\begin{rem}
By (\ref{e:blow1}) and (\ref{e:blow2}), we obtain that the
significant contributions  to integral $\int_\Omega f(u)dx$ come
from the largest core region and the boundary layers where $f$ is
large if $p=2$, but the core dominates for $p>2$.
\end{rem}
Let us consider two examples.

\noindent\textbf{Example 1.} Suppose $f(s)$ is decreasing,
$\int_0^\infty f(s)ds=1$, $f(s)\sim B/s^{1+b}$ as
$s\rightarrow\infty$ for some positive constants $b$ and $B$.

For $p=2$,
\[
  \frac{dM}{dt}\sim\frac{\Lambda_1^2}{f(M)},
\]
which implies
\[
  M\sim(\frac{b\Lambda_1^2}{B})^{-\frac{1}{b}}(T-t)^{-\frac{1}{b}}.
\]

For $p>2$,
\[
  \frac{dM}{dt}\sim\frac{\Lambda_2}{f^{p-1}(M)},
\]
which follows that
\[
  M\sim(\frac{(1+b)(p-1)-1}{B^{p-1}}\Lambda_2)^{\frac{1}{1-(1+b)(p-1)}}(T-t)^{\frac{1}{1-(1+b)(p-1)}}.
\]

\noindent\textbf{Example 2.} $f(s)=e^{-s}$.

For $p=2$,
\[
  \frac{dM}{dt}\sim\frac{\Lambda_1^2}{e^{-M}},
\]
which implies
\[
  M\sim-\ln(T-t)-2\ln\Lambda_1.
\]

For $p>2$,
\[
  \frac{dM}{dt}\sim\frac{\Lambda_2}{e^{(1-p)M}},
\]
that is,
\[
  M\sim
  \frac{1}{1-p}\ln\big((p-1)\Lambda_2\big)+\frac{1}{1-p}\ln(T-t).
\]


\begin{thebibliography}{99}
\bibitem{BL}J.W. Bebernes and A.A. Lacey,
            \textit{Global existence and finite-time blow-up for a class of nonlocal parabolic problems},
            Adv.  Differential Equations \textbf{2} (1997), 927-953.
\bibitem{C} J.A. Carrillo,
            \textit{On a nonlocal elliptic equation with decreasing nonlinearity
            arising in plasma physics and heat conduction},
            Nonlinear Anal. TMA \textbf{32} (1998), 97-115.
\bibitem{F} A.C. Fowler, I. Frigaard and S.D.  Howison,
            \textit{Temperature surges in current-limiting circuit devices},
             SIAM J.  Appl.  Math \textbf{52} (1992), 998-1011.
\bibitem{GT} D. Gilbarg and N.S. Trudinger,
            \textit{Elliptic Partial Differential Equations of Second Order},
            Springer-Verlag Berlin, 1977.
\bibitem{K} N.I. Kavallaris, A.A. Lacey and D.E. Tzanetis,
            \textit{Global existence and divergence of critical solutions
            of a nonlocal parabolic problem in Ohmic heating process},
             Nonlinear Anal.  TMA \textbf{58} (2004), 787-812.
\bibitem{KT} N.I. Kavallaris and D.E. Tzanetis,
            \textit{On the blow-up of the nonlocal thermistor problem},
            Proc.  Edinb.  Math.  Soc. \textbf{50} (2007), 389-409.
\bibitem{L1} A.A. Lacey,
            \textit{Thermal runaway in a nonlocal problem modelling Ohmic
            heating.  Part I: Model derivation and some special cases},
            European J.  Appl.  Math \textbf{6} (1995), 127-144.
\bibitem{L2} A.A. Lacey,
            \textit{Thermal runaway in a nonlocal problem modelling
            Ohmic heating.  Part II: General proof of blow-up and asymptotics of
            runaway},
            European J.  Appl.  Math \textbf{6} (1995), 201-224.
\bibitem{S} D.H. Sattinger,
            \textit{Monotone methods in nonlinear elliptic and parabolic boundary value problems},
            Indiana Univ.  Math.  J.  \textbf{21} (1972), 979-1000.
\bibitem{T} D.E. Tzanetis,
            \textit{Blow-up of radially symmetric solutions of nonlocal problem modelling Ohmic heating},
            Electron.  J.  Diff.  Eqns \textbf{11} (2002), 1-26.
\end{thebibliography}
 \end{document}